\documentclass[12pt,a4paper]{article}
\usepackage{amsmath}
\usepackage{amsthm}
\usepackage{amssymb}
\usepackage{amsbsy}
\usepackage{amsfonts}
\usepackage{amstext}
\usepackage{amscd}
\usepackage[dvips]{epsfig}

\numberwithin{equation}{section}
\theoremstyle{plain}
\newtheorem{thm}{Theorem}[section]
\newtheorem{prop}[thm]{Proposition}
\newtheorem{cor}[thm]{Corollary}

\theoremstyle{definition}

\newtheorem{rem}[thm]{Remark}
\newtheorem{defi}[thm]{Definition}

\DeclareMathOperator*{\real}{\mathbb{R}}

\DeclareMathOperator*{\im}{\text{Im}}

\newcommand{\pn}{\mathcal{P}(n)}
\newcommand{\ncpn}{\mathcal{NC}(n)}
\newcommand{\ipn}{\mathcal{I}(n)}
\newcommand{\mpn}{\mathcal{M}(n)}
\newcommand{\lpn}{\mathcal{LP}(n)}
\newcommand{\lncpn}{\mathcal{LNC}(n)}
\newcommand{\lipn}{\mathcal{LI}(n)}

\begin{document}
\title{The Monotone Cumulants}
\author{Takahiro Hasebe\footnote{Supported by Grant-in-Aid for JSPS Fellows. E-mail: hsb@kurims.kyoto-u.ac.jp} \\ and \\ Hayato Saigo\footnote{E-mail: harmonia@kurims.kyoto-u.ac.jp } \\ Graduate School of Science,  Kyoto University,\\  Kyoto 606-8502, Japan}
\date{}

\maketitle

\begin{abstract}
In the present paper we define the notion of generalized cumulants which gives a universal framework for 
commutative, free, Boolean, and especially, monotone probability theories. The uniqueness of generalized cumulants holds for each independence, and hence, 
generalized cumulants are equal to the usual cumulants in the commutative, free and Boolean cases. 
The way we define (generalized) cumulants needs neither partition lattices nor generating functions and then will give a new viewpoint to cumulants. We define ``monotone cumulants'' in the sense of generalized cumulants and we obtain quite simple proofs of  
central limit theorem and Poisson's law of small numbers in monotone probability theory. 
Moreover, we clarify a combinatorial structure of moment-cumulant formula with the use of ``monotone partitions''.  
\end{abstract}

keywords: Monotone independence; cumulants; umbral calculus 


\section{Introduction}
In quantum probability theory, many different notions of 
independence are defined. 
Among them, commutative, free, Boolean and monotone independence
 are considered as fundamental examples \cite{Mur3, Obata, S-W, V1}. 
 For commutative, free, Boolean and many other notions of independence, the 
associated cumulants and their appropriate 
generating functions have been introduced \cite{Leh1, Spe2, S-W, V2}. 
They are useful especially for the non-commutative versions of 
central limit theorems and Poisson's laws of small numbers. 

In the present paper we will introduce generalized cumulants which allow us to treat monotone independence. 
Since monotone independence depends on the order of random variables, 
the additivity of cumulants fails to hold. Instead, we introduce a generalized condition and then prove the uniqueness of generalized cumulants in Section \ref{gcumulants}. 
In Section \ref{mcumulants}, we show the existence of the monotone cumulants and obtain an 
explicit moment-cumulant formula for monotone independence. 
In Section \ref{limit}, we show the central limit theorem and Poisson's law of 
small numbers in terms of the monotone cumulants. 
  
The merit of the approach in this paper is that we do not need a partition lattice 
structure and generating functions to define cumulants. 
In addition, this approach to cumulants is also applicable to the theory of mixed cumulants 
for classical, free, Boolean and other notions of independence, which will be presented in another paper \cite{H-S}.  
We prove a combinatorial moment-cumulant formula in Section \ref{partition1}. 
This formula is expected to become a basis of combinatorics in monotone probability theory.  


\section{Four notions of independence}\label{independence}
Let $\mathcal{A}$ be a unital $*$-algebra over $\textrm{\boldmath $C$}$. A linear functional $\varphi: \mathcal{A} \to \textrm{\boldmath $C$}$ is called a state on $\mathcal{A}$ if 
$\varphi (a^* a) \geq 0$ and $\varphi (1_{\mathcal{A}})=1$. 
\begin{defi}
(1) An algebraic probability space is a pair $(\mathcal{A}, \varphi)$, where $\mathcal{A}$ is a unital $\ast$-algebra and $\varphi$ is a state on $\mathcal{A}$.  \\
(2) If the algebra $\mathcal{A}$ is a $C^\ast$-algebra, we call the pair $(\mathcal{A}, \varphi)$ a $C^\ast$-algebraic probability space. 
\end{defi}
An element $a$ in $\mathcal{A}$ is called an algebraic random variable. 
Quantities $\varphi (a_1a_2...a_n)$ are called mixed moments. 

The notion of independence in classical probability can be understood as a universal structure 
which gives a rule for calculating mixed moments, at least from the algebraic point of view. In quantum probability, lots of different notions of independence have been introduced. Among them, four notions mentioned below are known as fundamental examples \cite{Obata}. These four notions are important since they are ``universal products'' or ``natural products'' \cite{B-S1, Mur0, Mur4, Spe1}.

Let $(\mathcal{A},\varphi)$ be an algebraic probability space 
and $\left\{\mathcal{A}_{\lambda}; \lambda \in \Lambda  \right\}$ be a family of $*$-subalgebras of $\mathcal{A}$. In the following, 
 four notions of independence are defined as 
the rules for calculating mixed moments $\varphi (a_1a_2...a_n)$, where
\[
a_i\in \mathcal{A}_{\lambda _i},\:  a_i\notin \textrm{\boldmath $C$}1,\: \lambda_i\neq \lambda_{i+1},\: 1\leq i\leq n-1,\: n\geq 2. 
\]
\begin{defi}(Commutative independence).
$\left\{\mathcal{A}_{\lambda} \right\}$ is commutative independent if
\[
\varphi(a_1a_2...a_n)=\varphi(a_1)\varphi(a_2...a_n)
\]
holds when $\lambda_1\neq \lambda_r$ for all $2\leq r\leq n$, and otherwise, 
letting $r$ be the least number such that 
$\lambda_1= \lambda_r$,
\[
\varphi(a_1a_2...a_n)=\varphi(a_2...a_{r-1}(a_1a_r)a_{r+1}...a_n). 
\]  
\end{defi}
\begin{defi}(Free independence \cite{V1}).
$\left\{\mathcal{A}_{\lambda} \right\}$ is free independent if
\[
\varphi(a_1a_2...a_n) = 0 
\]
holds whenever $\varphi (a_1)=...=\varphi (a_n)=0$.
\end{defi}
\begin{defi}(Boolean independence \cite{S-W}).
$\left\{\mathcal{A}_{\lambda} \right\}$ is Boolean independent if
\[
\varphi(a_1a_2...a_n)=\varphi(a_1)\varphi(a_2...a_n). 
\]
\end{defi}
\begin{defi}(Monotone independence \cite{Mur3}).
Assume that the index set $\Lambda $ is equipped with a linear order $<$.
$\left\{\mathcal{A}_{\lambda} \right\}$ is monotone independent if
\[
\varphi (a_1...a_i...a_n)=\varphi (a_i)\varphi(a_1...a_{i-1}a_{i+1}...a_n)
\]
holds when $i$ satisfies $\lambda _{i-1}<\lambda _i$ and $\lambda _i>\lambda _{i+1}$ (one of the inequalities is eliminated when $i=1$ or $i=n$).
\end{defi}
 
A system of algebraic random variables $\left\{x_\lambda  \right\}$ are called commutative/free/Boolean/ \\ 
monotone independent if $\left\{\mathcal{A}_\lambda \right\}$ are 
commutative/free/Boolean/monotone independent, where $\mathcal{A}_\lambda$ denotes the algebra generated by $x_\lambda$.

\section{Generalized cumulants}\label{gcumulants}
Let $(\mathcal{A}, \varphi)$ be an algebraic probability space. 
We define $M_n(X):= \varphi(X^n)$.  
Important properties of cumulants $\{K_n \}_{n \geq 1}$ for any one of commutative, free and Boolean independence are summarized as follows \cite{Leh1}: 
\begin{itemize}
\item[(K1)] Additivity:  If $X$,  $Y \in \mathcal{A}$ are (commutative, free or Boolean) independent,  
\begin{equation}\label{addicum}
K_n(X + Y) = K_n(X) + K_n(Y). 
\end{equation}
for any $n \geq 1$. \\
\item[(K2)] Homogeneity: for any $\lambda > 0$ and any $n$,  
\begin{equation}\label{1a}
K_n(\lambda X) = \lambda ^n K_n(X). 
\end{equation} 
\item[(K3)] For any $n$, there exists a polynomial $Q_n$ of $n-1$ variables such that 
\begin{equation}\label{1b}
M_n(X) = K_n(X) + Q_n(K_1(X), \cdots, K_{n-1}(X)), 
\end{equation}
where $M_n(X)$ is the $n$-th moment of $X$. 
\end{itemize}
We introduce generalized cumulants which allow us to treat monotone independence. Since monotone independence depends on the order of random variables, the additivity of cumulants fails to hold. Instead, we introduce a generalized condition: 
\begin{itemize}
\item[(K1')] If $X^{(1)}, X^{(2)}, \cdots, X^{(N)}$ are independent, 
identically distributed to $X$ in the sense of moments (i.e., $M_n(X)=M_n(X^{(i)})$ for all $n$),  we have
\begin{equation}\label{1c}
K_n(N.X):= K_n(X^{(1)} + \cdots + X^{(N)}) = N K_n(X).
\end{equation}
\end{itemize}
Here we understand that $K_n(0.X):=\delta_{n0}$. We note that the notation $N.X$ 
is inspired by the ``dot product'' or ``random sum'' operation in 
the theory of the classical umbral calculus \cite{DiNardo, Rota}. (Probably the notion above 
will be used as a foundation for ``the non-commutative umbral calculus''.)  

We show the uniqueness of generalized cumulants w.r.t. the notion of each independence. 
\begin{thm}
Generalized cumulants satisfying (K1'), (K2) and (K3) are unique and the $n$-th cumulant is given by 
the coefficient of $N$ in $M_n(N.X)$. 
\end{thm}
\begin{proof}
By (K3) and (K1), we obtain 
\begin{equation}\label{cum4}
\begin{split}
M_n(N.X) &= K_n(N.X) +Q_n (K_1(N.X), \cdots, M_{n-1}(N.X)) \\
         &= NK_n(X) + Q_n(NK_1(X), \cdots, NK_{n-1}(X)). 
\end{split}
\end{equation} 
Therefore, $M_n(N.X)$ is a polynomial of $N$ and $M_k(X)$ $(1 \leq k \leq n)$. 
By condition (K2), the polynomial $Q_n$ does not contain linear terms and a constant for any $n$; hence the coefficient of the linear term $N$ is nothing but $K_n(X)$. Uniqueness of cumulants also follows from the above observation. More precisely, if another cumulants $K_n '$ are given, 
there exist polynomials $Q_n '$ in the condition (K3) for $K_n '$. Then we have
 \begin{equation*}
\begin{split}
M_n(N.X) &=  NK_n(X) + Q_n(NK_1(X), \cdots, NK_{n-1}(X)) \\
               &=  NK_n ' (X) + Q_n ' (NK_1 '(X), \cdots, NK_{n-1} '(X)).   
\end{split}
\end{equation*} 
This is an identity of polynomials of $N$; therefore, the coefficients of $N$ coincide and $K_n(X) = K_n '(X)$ holds.   
\end{proof}
From now on, we use the word ``cumulants'' instead of ``generalized cumulants''  to label $K_n$ above. 

\section{The monotone cumulants}\label{mcumulants}
The Cauchy transformation of a random variable $X$ is defined by 
\[
G_X(z):=\sum_{n = 0}^\infty \frac{M_n(X)}{z^{n+1}}. 
\]
We consider this in the sense of a formal power series if the series is not absolutely  convergent. Muraki proved in \cite{Mur3} that $G_{X+Y}(z) = G_X(\frac{1}{G_Y(z)})$ holds if $X$ and $Y$ are  monotone independent. We give a simple derivation of the formula.  
\begin{prop}
For monotone independent random variables $X$ and $Y$, it holds that 
\begin{equation}\label{moment6}
\begin{split}
M_n(X + Y) &= \sum_{k = 0} ^{n} \sum_{\substack{j_0 + j_1 + \cdots +j_k = n - k, \\  0 \leq j_l, ~0 \leq l \leq k }} M_k(X) M_{j_0}(Y)\cdots M_{j_k}(Y) \\ 
           &= M_n(X) + M_n(Y) + \sum_{k = 1} ^{n-1} \sum_{\substack{j_0 + j_1 + \cdots +j_k = n - k, \\  0 \leq j_l, ~0 \leq l \leq k }} M_k(X) M_{j_0}(Y)\cdots M_{j_k}(Y). 
\end{split}
\end{equation} 
\end{prop}
\begin{proof}
$(X + Y)^n$ can be expanded as 
\begin{equation}
\begin{split}
(X + Y)^n &= X^n + Y^n + \sum_{k = 1} ^{n-1} \sum_{\substack{j_0 + j_1 + \cdots +j_k = n - k, \\  0 \leq j_l, ~0 \leq l \leq k }} Y^{j_0}XY^{j_1}X \cdots X Y^{j_k}. 
\end{split}  
\end{equation}
Taking the expectation of the above equality, we obtain (\ref{moment6}). 
\end{proof}
\begin{cor}\label{cor11}
There exists a polynomial $P_n ^M$ of $2n-2$ variables for any $n \geq 1$ such that 
\begin{equation}\label{moments12}
M_n(X + Y) = M_n(X) + M_n(Y) + P_n ^M(M_1(X), \cdots, M_{n-1}(X),M_1(Y), \cdots, M_{n-1}(Y)) 
\end{equation}
holds if $X$ and $Y$ are monotone independent.  
\end{cor}
\begin{rem}\label{poly0}
A similar result is valid for any other independence: there exists a polynomial $P_n ^C$ (resp. $P_n^F, P_n ^B$) such that (\ref{moments12}) holds, with $P_n ^M$ replaced by another polynomial $P_n ^C$ (resp. $P_n ^F, P_n ^B$),  if $X$ and $Y$ are commutative (resp. free and Boolean) independent.
\end{rem}
By this corollary, we obtain the proposition below.
\begin{prop}\label{poly1}
$M_n(N.X)$ is a polynomial of $N$ (without a constant term) for any $n\geq 0$.
\end{prop}
\begin{proof}
We use induction w.r.t. $n$. For $n=1$, it is obvious from linearity of expectation.
Suppose the proposition holds for $n\leq l$. From (\ref{moments12}), we obtain 
\begin{equation*}
\Delta M_{l+1}(N.X)= M_{l+1}(X)+ P_n ^M(M_1(X), \cdots, M_{l}(X), M_1((N-1).X),\cdots, M_{l}((N-1).X)). 
\end{equation*}
Here, $\Delta M_{l+1}(N.X):= M_{l+1}(N.X)-M_{l+1}((N-1).X)$. By the assumption of induction, 
$P_{l+1} ^M(M_1(X), \cdots, M_{l}(X), M_1((N-1).X),\cdots, M_{l}((N-1).X))$ is a polynomial of $N$. Then $M_{l+1}(N.X)$ is a polynomial of $N$ 
(without a constant term) 
because $\Delta M_{l+1}(N.X)$ is a polynomial of $N$ and $M_{l+1}(0.X)=0$. 
\end{proof}
 
As the proposition above holds, we may define $m_n(t)=M_n(t.X)$ by replacing $N$ with $t \in \real$.
Note that this is a polynomial w.r.t. $t$ and that $m_n(1)=M_n(X)$. Moreover, we easily obtain 
\begin{equation}\label{flow}
m_n(t+s) = 
m_n(t) + m_n(s) + 
\sum_{k = 1} ^{n-1} 
\sum_{\substack{j_0 + j_1 + \cdots +j_k = n - k, 
\\  0 \leq j_l, ~0 \leq l \leq k }} m_k(t) m_{j_0}(s)\cdots m_{j_k}(s) 
\end{equation}
from the definition of $m_n(t)$ and (\ref{moment6}).

Now we come to define the main notion.
\begin{defi}\label{defi11}
Let $r_n=r_n(X)$ be the coefficient of $N$ in $M_n(N.X)$ (or the coefficient of $t$ in $m_n(t)$). 
We call $r_n$ the $n$-th monotone cumulant of $X$.
\end{defi}
\begin{rem}
(1) A result analogous to Corollary \ref{cor11} is also valid for any one of commutative, free and Boolean independence as mentioned in Remark \ref{poly0}, and therefore, we can prove Proposition \ref{poly1} for any independence. This enables us to define cumulants as the coefficient of $N$ in $M_n(N.X)$ also for commutative, free and Boolean independence. This is a simple definition of cumulants 
for each independence without use of generating functions.   \\
(2) There is an interesting formula called Good's formula in commutative, free and Boolean cases \cite{Leh1}.  Lehner defined mixed cumulants 
for the three notions of independence in terms of Good's formula.  The monotone independence is, however,  non-commutative and Lehner's approach cannot be directly applied to monotone cumulants.  
\end{rem} 
The monotone cumulants satisfy the axioms (K1') and (K2) 
because $M_n(N.(M.X))=M_n((NM).X)$ and $M_n(N.(\lambda X))=M_n(\lambda (N.X))$. The former equality is a consequence of the associativity of monotone independence or monotone convolution. 
We prepare the following proposition for the moment-cumulant formula. 

\begin{prop}
The equations below hold: 
\begin{equation}\label{system of ODEs}
\begin{split}
   \frac{dm_0(t)}{dt} &= 0, \\
\frac{dm_n(t)}{dt} &= \sum_{k=1} ^n k r_{n-k + 1} m_{k-1}(t) \text{~~for $n \geq 1$}, 
\end{split}
\end{equation}
with initial conditions $m_0(0) = 1$ and $m_n (0) = 0$ for $n \geq 1$. 
\end{prop}
\begin{proof}
From (\ref{flow}), we obtain
\begin{equation}\label{flow to derivative}
m_n(t+s)- 
m_n(t) = m_n(s) + 
\sum_{k = 1} ^{n-1} 
\sum_{\substack{j_0 + j_1 + \cdots +j_k = n - k, 
\\  0 \leq j_l, ~0 \leq l \leq k }} m_k(t) m_{j_0}(s)\cdots m_{j_k}(s). 
\end{equation}
By definition,
\begin{equation}
m_i(s)=r_is+s^2(\cdots )
\end{equation}
holds. Comparing the coefficients of $s$ in (\ref{flow to derivative}), we obtain the conclusion. 
\end{proof}

We show that $\{M_n(X)\}_{n \geq 0}$ and $\{r_n(X) \}_{n \geq 1}$ are connected with each other by a formula. 
\begin{thm}\label{m-c121}  The following formula holds:   
\begin{equation}\label{m-c122}
 M_n(X) = \sum_{k =1} ^{n} \sum_{1 = i_0 < i_1 < \cdots < i_{k-1} < i_k = n +1} 
\frac{1}{k!}\prod_{l = 1} ^{k} i_{l-1} r_{i_l - i_{l-1}}(X). 
\end{equation} 
\end{thm}
\begin{proof}
This formula is obtained directly by (\ref{system of ODEs}). 
We shall use the equations in the integrated forms 
\begin{equation}\label{integrated}
\begin{split}
m_0(t) &= 1, \\
m_n(t)   &= \sum_{k=1} ^n k r_{n-k + 1} \int_{0} ^{t} m_{k-1}(s) ds \text{~~for $n \geq 1$}.
\end{split}
\end{equation}
Then we have 
\begin{equation*}
\begin{split}
m_n (t) &= \sum_{k_1 =1} ^n k_n r_{n-k_1 + 1} \int_{0} ^{t} m_{k_1 -1}(t_1) dt_1  \\
        &= \sum_{k_1 = 1} ^n \sum_{k_{2} = 1} ^{k_1 -1}k_1 k_2 r_{n-k_1 +1} r_{k_1 -k_2}  \int_0 ^{t}dt_1 \int_0 ^{t_1} dt_2 ~m_{k_2 -1}(t_2) \\ 
        &= \sum_{k_1 = 1} ^n \sum_{k_2 = 1} ^{k_1 -1} \sum_{k_3 = 1} ^{k_2 -1} 
k_1 k_2 k_3 r_{n-k_1 +1} r_{k_1 -k_2} r_{k_2 - k_3}   \int_0 ^t dt_1 \int_0 ^{t_1} dt_2 \int_0 ^{t_2}dt_3 ~m_{k_3 -1}(t_3)dt_3 \\ 
       &= \cdots. \\
\end{split}
\end{equation*}
When this calculation ends, we obtain the formula 
\[
m_n(t) = \sum_{k =1} ^{n} \sum_{1 = i_0 < i_1 < \cdots < i_{k-1} < i_k = n +1} 
\frac{t^k}{k!}\prod_{l = 1} ^{k} i_{l-1} r_{i_l - i_{l-1}}, 
\]
where $i_l:=k_{n-l}$. Putting $t = 1$, we have (\ref{m-c122}).  \\
\end{proof}
\begin{rem}
This formula has been already obtained in the case of the monotone Poisson distribution \cite{Belt1, Belt2}. 
\end{rem}
\begin{cor}
The monotone cumulants $r_n=r_n(X)$ satisfy (K3).    
\end{cor}

Hence, we obtain the main theorem.
\begin{thm}
$r_n$ are the unique (generalized) cumulants for monotone independence.
\end{thm}

\section{Limit theorems in monotone probability theory}\label{limit}
As applications of monotone cumulants, we give short proofs of limit theorems which have been obtained by combinatorial arguments in \cite{Mur2}. 
\begin{thm}\label{central} Let $(\mathcal{A}, \phi)$ be a $C^\ast$-algebraic probability space. 
Let $X^{(1)}, \cdots, X^{(N)}, \cdots$ be identically distributed, monotone independent self-adjoint random variables with $\phi(X^{(1)}) = 0$ and $\phi((X^{(1)})^2) = 1$. Then the probability distribution of $X_N:= \frac{X^{(1)} + \cdots + X^{(N)} }{\sqrt{N}}$ converges weakly to the arcsine law with mean $0$ and variance $1$. 
\end{thm} 
\begin{proof}
By the properties (K1') and (K2), we immediately obtain $r_1(X_N)=0$, $r_2(X_N) = 1$ and $r_n(X_N) = N^{-\frac{n-2}{2}}r_n(X^{(1)}) \to 0$ as $N \to \infty$. 
By (K3), $M_n(X_N)$ converges to $M_n$ characterized by the monotone cumulants $(r_1, r_2, r_3, r_4, \cdots) = (0, 1, 0, 0, \cdots)$. 
Since only $r_2 (=1)$ is nonzero in (\ref{m-c122}), we can calculate the moments as $M_{2n-1} = 0$ and $M_{2n} = \frac{(2n-1)!!}{n!}$ for all $n \geq 1$, where the double factorial $(2n-1)!!$ is defined by $1\cdot 3  \cdots (2n-3)(2n-1)$ for $n \geq 1$.  
The limit measure is the arcsine law with mean $0$ and variance $1$: 
\[
\int_{-\sqrt{2}}^{\sqrt{2}}\frac{x^{2n}}{\pi\sqrt{2-x^2}}dx = \frac{(2n-1)!!}{n!}. 
\] The moment problem of the arcsine law is determinate and therefore the distribution of $X_N$ converges to the arcsine law weakly (see Theorem 4.5.5 in \cite{KLC}).  
\end{proof}
We can show Poisson's law of small numbers similarly in the setting of a triangular array.
\begin{thm}\label{poisson}
Let $X_N ^{(n)}$ $(1 \leq n \leq N,~1 \leq N < \infty)$ be self-adjoint random variables in a $C^\ast$-algebraic probability space such that \\ 
(1) for each $N$, $X_N ^{(n)} (1 \leq n \leq N)$ are identically distributed, monotone independent self-adjoint random variables;  \\ 
(2) $N M_k(X_N^{(1)}) \to \lambda > 0$ as $N \to \infty$ for all $k \geq 1$. \\ Then the distribution of $X_N := X_N^{(1)} + \cdots + X_N^{(N)}$ converges weakly to the monotone Poisson distribution with parameter $\lambda$. 
\end{thm}
\begin{proof}
By properties (K1') and (K3), $r_n(X_N) = Nr_n(X_N^{(1)}) = NM_n(X_N^{(1)}) + o(1) \to \lambda$ for $n \geq 1$. Here we used the fact that the polynomial in (K3) does not contain linear terms of $M_k(X)$ $(1 \leq k \leq n-1)$. Therefore, the limit moment $M_n$ is characterized by the monotone cumulants $(r_1, r_2, r_3, \cdots) = (\lambda, \lambda, \lambda, \cdots)$. From (\ref{m-c122}) we have 
\[
M_n = \sum_{k =1} ^{n}\frac{\lambda^k}{k!} \sum_{1 = i_0 < i_1 < \cdots < i_{k-1} < i_k = n +1} i_0 i_1 \cdots i_{k-1}.  
\]
It is known that this gives a determinate moment sequence and the limit distribution is called the monotone Poisson distribution (see \cite{Belt1, Belt2, Mur3, Mur2}). Therefore, the distribution of $X_N$ converges weakly to the monotone Poisson distribution. 
\end{proof}

If we formulate the the above theorems in terms of monotone convolution $\rhd$ of  probability measures, we can include probability measures with possibly noncompact supports. We now explain this. 
\begin{defi} \cite{Mur3}
The monotone convolution $\mu \triangleright \nu$ of probability measures $\mu$ and $\nu$ is define by the relation 
\[
H_{\mu \triangleright \nu} (z) = H_\mu \circ H_\nu (z),~~~\im ~z \neq 0,  
\]
where $H_\mu$ is defined by $H_\mu(z) = [\int_{\real}\frac{\mu(dx)}{z-x}]^{-1}$. $H_\mu$ is called the reciprocal Cauchy transform of $\mu$.  
\end{defi}
Then the definition of cumulants $r_n(\mu)$ of a probability measure $\mu$ with finite moments is basically the same as that of Definition \ref{defi11}; the $n$-th monotone cumulant $r_n(\mu)$ is defined as the coefficient of $N$ in $m_n(\mu^{\rhd N})$. 

The dilation operator $D_\lambda$ is defined so that $\int_{\real} f(x) D_\lambda \mu(dx) = \int_{\real} f(\lambda x) \mu(dx)$ holds for all bounded continuous function $f$ on $\real$. 
\begin{thm} \label{central2}
(1) Let $\mu$ be a probability measure on $\real$ with finite moments of all orders, mean $0$ and variance $1$. 
Then the probability measure $D_{\frac{1}{\sqrt{N}}}\mu \triangleright \cdots \triangleright D_{\frac{1}{\sqrt{N}}}\mu$ ($N$ times) converges weakly to the arcsine law with mean $0$ and variance $1$. \\
(2) Let $\mu^{(N)}$ be probability measures on $\real$ with finite moments of all orders. 
We assume that  $NM_k(\mu^{(N)}) \to \lambda > 0$ as $N \to \infty$ for every $k \geq 1$. Then $\mu_N := \mu^{(N)} \triangleright \cdots \triangleright \mu^{(N)}$ ($N$ times) converges to the monotone Poisson distribution with parameter $\lambda$. 
\end{thm}

The proof is totally identical to those of Theorem \ref{central} and Theorem \ref{poisson} if we replace respectively the moments  and monotone cumulants of random variables with those of probability measures. 
This is because the convergence of moments implies 
the weak convergence of probability measures, if the limit moments are determinate; one need not assume that the initial 
sequences of the probability measures $D_{\frac{1}{\sqrt{N}}}\mu$ and $\mu^{(N)}$ have determinate 
moments (see Theorem 4.5.5 in \cite{KLC}). 
Therefore, the proof of Theorem \ref{central2} is identical to those of Theorem \ref{central} and Theorem \ref{poisson}.

\section{Moment-cumulant formula by monotone partitions}\label{partition1}
In the classical, free and Boolean cases, the moment-cumulant formulae are described with the use of 
the structures of partitions. Let $\pn$ be the set of all partitions of $\{1, \cdots, n \}$ \cite{Shi}, let 
$\ncpn$ be the set of all non-crossing partitions of $\{1, \cdots, n \}$ \cite{Spe2} and let $\ipn$ 
be the set of all interval partitions \cite{S-W}. We denote the number of the elements in a set $V$ by $|V|$. For a sequence of real numbers $\{t_n \}_{n \geq 1}$ and $\pi = \{V_1, \cdots, V_l\} \in \pn$ we define $t(\pi) := \prod_{j = 1}^l t_{|V_j|}$. Then the moment-cumulant formulae are written in the forms
\begin{align}
&m_n = \sum_{\pi \in \pn} r(\pi) \text{~~~~(classical case)}, \label{cmc} \\
&m_n = \sum_{\pi \in \ncpn} r(\pi) \text{~~~~(free case)}, \label{fmc} \\
&m_n = \sum_{\pi \in \ipn} r(\pi) \text{~~~~(Boolean case)}. \label{bmc} 
\end{align}
These sets of partitions appear also as the highest coefficients of the decomposition rules of "universal products" in the classification of independence \cite{B-S1, Spe1}. Connections between the formulae (\ref{cmc})-(\ref{bmc}) and 
the highest coefficients seem to be not known yet. 
\begin{figure}
\centering
\includegraphics[width=8cm,clip]{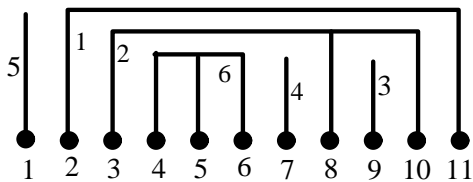}
\caption{An example of a monotone partition.}
\label{dia1}
\end{figure}

Muraki has defined the notion of linearly ordered partitions and classified quasi-universal products. 
Let $(\pi, \lambda)$ be a pair which consists of $\pi \in \pn$ and a linear ordering $\lambda$ of the blocks of $\pi$. 
It is useful to denote $(\pi, \lambda)$ by $\pi = \{V_1 < V_2 < \cdots < V_l\}$. He has introduced the sets 
$\lpn$, $\lncpn$ and $\lipn$, where $\mathcal{L}$ denotes the structure of linear orderings. 
We note that $|\lpn| = \sum_{k = 1} ^{n} k!|\{\pi \in \pn; |\pi| = k \}|$, for instance. A linearly ordered partition can be visualized by a diagram with blocks labeled by natural numbers. For instance, Fig. \ref{dia1} describes the partition 
$\{\{2, 11 \} < \{3, 8, 10 \} < \{9 \} < \{7 \} < \{1 \} < \{4, 5, 6 \}\}$. 
He has defined the set of monotone partitions by 
\begin{equation}
\mpn:= \{(\pi, \lambda); \pi \in \ncpn,~ \text{if $V, W \in \pi$ and $V$ is in the inner side of $W$, then $V >_{\lambda} W$} \} 
\end{equation}
which appears as the highest coefficients of the decomposition rules of monotone product. Fig. 1 is an example of a  monotone partition. 

Later Sa\l apata and Lenczewski have introduced the same partitions from a different viewpoint. They have defined the notion of $m$-monotone in \cite{Len1} as an interpolation between monotone and free probabilities, and used a generalization of monotone partitions to compute the moments of limit distributions of central limit theorem and Poisson's law of small numbers. The $m = 0$ case corresponds to the monotone partitions. 

We have discovered that the monotone partitions play crucial role in the context of monotone cumulants as follows.   
\begin{thm}
Let $r_n$ be the monotone cumulants of a probability measure with finite moments of all orders. 
Then the formula  
\begin{equation}
m_n = \sum_{(\pi, \lambda) \in \mpn} \frac{r(\pi)}{|\pi|!} \label{mmc}  
\end{equation}
holds. 
\end{thm} 
Before proving this, we note that this formula can be understood naturally since the formulae (\ref{cmc})-(\ref{bmc}) are rewritten in the forms 
\begin{align}
&m_n = \sum_{(\pi, \lambda) \in \lpn} \frac{r(\pi)}{|\pi|!} \text{~~~~(classical case)}, \label{cmc2}\\
&m_n = \sum_{(\pi, \lambda) \in \lncpn} \frac{r(\pi)}{|\pi|!} \text{~~~~(free case)}, \label{fmc2}\\
&m_n = \sum_{(\pi, \lambda) \in \lipn} \frac{r(\pi)}{|\pi|!} \text{~~~~(Boolean case)}. \label{bmc2} 
\end{align}

\begin{figure}
\centering
\includegraphics[width=3cm,clip]{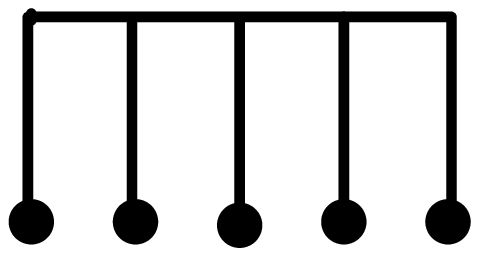}~~~~~~~~~~~~~~~~~
\includegraphics[width=2cm,clip]{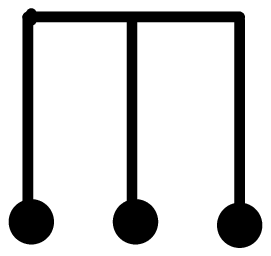}
\caption{Two examples of blocks of interval type.}
\label{inter1}
\end{figure}
\begin{proof}
Let $\pi \in \pn$. 
We call a block $V \in \pi$ a block of \textit{interval type} if $V$ is of such a form as $V = \{j, j+1, \cdots, j + k \}$ for $1 \leq j \leq n$, $k \geq 0$ (see Fig. \ref{inter1}). In other words, $V$ is a block of interval type if $V$ does not contain other blocks in the inner side of itself. 
Let $k \geq 1$ be an integer. For a given monotone partition $\pi = \{V_1 < \cdots < V_k \} \in \mpn$ (with $|\pi| = k$) we can define a map $T_k: \{V_1 < \cdots < V_k \} \mapsto (|V_1|, \cdots, |V_k|)$. 
Conversely, for a given sequence of integers $(n_1, \cdots, n_k)$ with $n_1 + \cdots + n_k = n$ and $n_j \geq 1$, we need to count the number $|T_k ^{-1}(n_1, \cdots, n_k)|$ to prove the moment-cumulant formula.  
We now consider a procedure for producing all monotone partitions in $T_k ^{-1}(n_1, \cdots, n_k)$. 
Before the definition of the procedure, we note important properties of monotone partitions. 
Let $\{V_1 < \cdots < V_k \}$ be a monotone partition. First we note that $V_k$ is a block of interval type since 
$V_k$ does not contain other blocks in the inner side of itself. Second, for any $1 \leq  j \leq k$, $V_j$ can be seen as a block of interval type if we forget the higher blocks $V_{j+1}, \cdots, V_{k}$. For instance, Fig. \ref{dia3} is a diagram in Fig. \ref{dia1} without the blocks $\{9 \}, \{7 \}, \{1 \}, \{4, 5, 6 \}$. Then the block $\{3, 8, 10 \}$ can be seen as a block 
of interval type in Fig. \ref{dia3}. 
\begin{figure}
\centering
\includegraphics[width=8cm,clip]{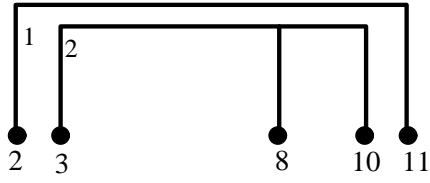}
\caption{A diagram in Fig. \ref{dia1} without the blocks $\{9 \}, \{7 \}, \{1 \}, \{4, 5, 6 \}$.}
\label{dia3}
\end{figure}
Taking this property into consideration, we consider the following procedure. 
The key point is to choose the blocks in order from the highest one to the lowest one. 
\begin{itemize}
\item[(1)] Choose $1 \leq k \leq n$ and $(n_1, \cdots, n_k)$ with $n_1 + \cdots + n_k = n$ and $n_j \geq 1$ and fix them.  
\item[(2)] Choose the position of the block $V_k$ of interval type (with $|V_k| = n_k$) among the $n$ elements and remove the block $V_k$. Then there remain $n - n_k$ elements. 
\item[(3)] Repeat the procedure (2) similarly for the block $V_{k-1}$ of interval type. Then the remaining elements are 
$n - n_k - n_{k-1}$. 
\item[(4)] Similarly, we choose the blocks $V_{k-2}, \cdots, V_1$.  These blocks, equipped with the linear ordering $V_1 < \cdots < V_k$, determine a monotone partition $(\pi, \lambda)$.  
\end{itemize}

\begin{figure}
\centering
\includegraphics[width=4.5cm,clip]{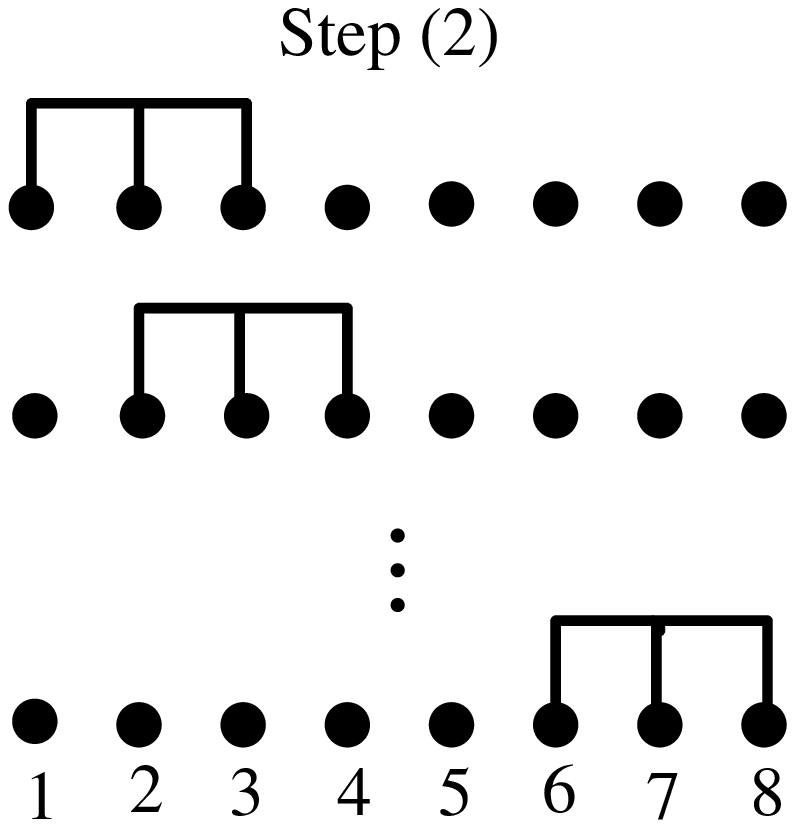}
\includegraphics[width=3cm,clip]{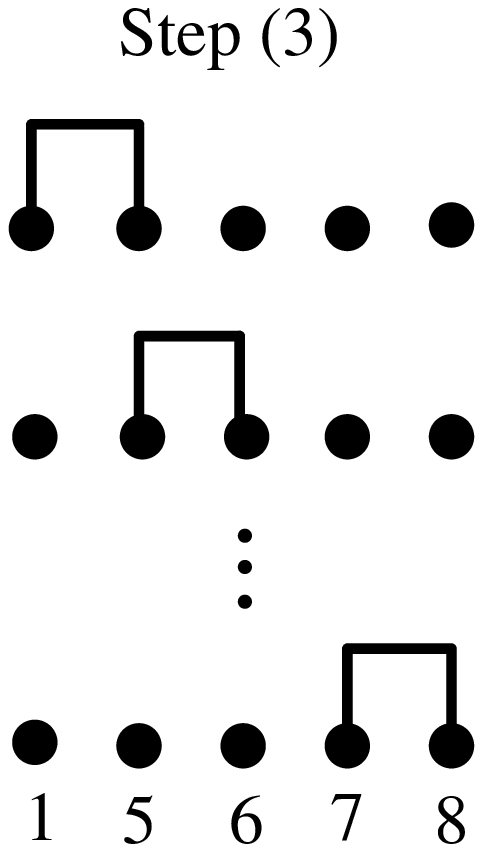}
\caption{An example of the steps (2) and (3).}
\label{dia4}
\end{figure}
We explain the above through an example. We put $n = 8$, $k \geq 3$, $n_k = 3$ and $n_{k-1} =2$. 
In step (2), there are six possibilities to choose $V_k$ in Fig. \ref{dia4}. For instance, when we choose $V_k =\{2, 3,4 \}$, after the removal of $V_k$ we choose $V_{k-1}$ in Step (3). Then there are four possibilities to choose 
$V_{k-1}$.   

It is not difficult to see that every linearly ordered partition in $T_k ^{-1}(n_1, \cdots, n_k)$ appears just once in the above procedure (2)-(4). Moreover, we can count the number as follows. In Step (2) there are $n - n_k + 1$ ways to choose the position of $V_k$. In Step (3) there are $n-n_k-n_{k-1} +1$ ways to choose the position of $V_k$. Similarly, we can count the number of all possible ways to choose the blocks $V_k, \cdots, V_1$, which is equal to $(n-n_k + 1)(n-n_k - n_{k-1} +1) \cdots (n - \sum_{j = 1}^{k} n_j +1)$. 
Therefore, we have 
\begin{equation}
\begin{split}
\sum_{(\pi, \lambda) \in \mpn} \frac{r(\pi)}{|\pi|!} 
        &= \sum_{k = 1}^n \sum_{\substack{n_1 + \cdots + n_k = n, \\n_j \geq 1,~ 1 \leq j \leq k}} \frac{1}{k!}\Big{(}\prod_{m = 1}^{k}(n - \sum_{j = m}^{k} n_j +1)\Big{)} r_{n_1} \cdots r_{n_k} \\
        &=  \sum_{k = 1}^n \sum_{1 = i_0 < i_1 < \cdots < i_{k-1} < i_{k}= n+1} \frac{i_1 i_2 \cdots i_{k-1}}{k!} r_{i_1 - i_0} \cdots r_{i_k - i_{k-1}},  
\end{split}
\end{equation}
where we have put $i_m := n - \sum_{j = m+1}^{k}n_j + 1$. The last expression is equal to $m_n$ by Theorem \ref{m-c121}. 
\end{proof}
\begin{rem}
The moment-cumulant formula includes the cases of the normalized arcsine law $(r_1, r_2, r_3, r_4, \cdots) = (0,1,0,0,\cdots)$ and the monotone Poisson distribution $(r_1, r_2, r_3, r_4, \cdots) = (\lambda,\lambda,\lambda,\lambda,\cdots)$ obtained in \cite{Len1, Len2} as $m = 0$ case. 
See Theorem 7.1 and Theorem 8.1 in \cite{Len1}.
The prototype of the combinatorial discussion in the proof above is in \cite{Saigo}.  
\end{rem}

\section*{Acknowledgements} 
The authors are grateful to Prof. 
Nobuaki Obata for guiding them to an important reference \cite{KLC} and encouragements. 
 They thank the referee for many constructive comments. 
 They also thank Prof. Izumi Ojima, Prof. Shogo Tanimura,
  Mr. Ryo Harada, Mr. Hiroshi Ando, Mr. Kazuya Okamura and Mr. Daisuke Tarama for their interests and comments. 
They are grateful to Prof. Rafa\l \ Sa\l apata for indicating the works \cite{Len1, Len2}. This work was supported by JSPS KAKENHI 21-5106.

\end{document}